\documentclass[12pt]{amsart}
\usepackage{amsfonts, amssymb,amsmath, amsthm, amsxtra, latexsym, amscd}
\usepackage[latin1]{inputenc} 

\def\draft{\centerline{(Draft {\the \day}/{\the\month} \the \year.)}}

\theoremstyle{definition}
\newtheorem{theo+}    {Theorem}      [section]
\newtheorem{prop+}  [theo+]  {Proposition}
\newtheorem{coro+}  [theo+]  {Corollary}
\newtheorem{lemm+}  [theo+]  {Lemma}
\newtheorem{deep+}  [theo+]  {Deep Result}
\newtheorem{fact+}  [theo+]  {Fact}
\theoremstyle{definition}
\newtheorem{exam+}  [theo+]  {Example}
\newtheorem{rema+}  [theo+]  {Remark}
\newtheorem{defi+}  [theo+]  {Definition}
\newtheorem{xca+}[theo+]{Exercise}

\numberwithin{equation}{section}

\def\draft{\centerline{(Draft {\the \day}/{\the\month} \the \year.)}}
\bigskip

\def\refn#1.#2{\expandafter\def\csname#1\endcsname{[#2]}}
\def\refnr#1.{\csname#1\endcsname}

\def\a{\alpha}

\def\Claminv2{|C(\Lambda)|^{-2}}

\def\varepsi{\varepsilon}

\def\Ome{\Omega}

\def\de{d\varepsilon}

\def\Aa2D{A^{\a,2}(D)}
\def\bAa2D{\overline{A^{\a,2}(D)}}
\def\Ab2D{A^{\beta,2}(D)}
\def\bAb2D{\overline{A^{\beta,2}(D)}}

\def\Norm#1_#2{\Vert#1\Vert_{#2}}

\def\phipl12{\phi_{p_{l_1}, p_{l_2}}}
\def\phip01{\phi_{p_{0}, p_{0}}}

\def\a{\alpha}

\def\Claminv2{|C(\Lambda)|^{-2}}

\def\varepsi{\varepsilon}

\def\sig{\sigma}

\def\Ker{\operatorname{Ker}}

\def\de{d\varepsilon}

\def\Aa2D{A^{\a,2}(D)}
\def\bAa2D{\overline{A^{\a,2}(D)}}
\def\Ab2D{A^{\beta,2}(D)}
\def\bAb2D{\overline{A^{\beta,2}(D)}}

\def\phipl12{\phi_{p_{l_1}, p_{l_2}}}
\def\phip01{\phi_{p_{0}, p_{0}}}

\def\alg/{algebra} 
\def\Alg/{Algebra} 
\def\alt/{alternative} 
\def\anal/{analytic}
\def\analfunc/{\anal/\ \func/}
\def\Ans/{\it Answer. \normal}
\def\ass/{associative}
\def\nass/{non-\ass/}
\def\autom/{automorphism}
\def\homom/{homomorphism}
\def\isom/{isomorphism}
\def\bdd/{bounded}
\def\Bdd/{Bounded}
\def\bddsymdom/{bounded \sym/ \dom/}
\def\Cartdom/{Cartan \dom/}
\def\bdry/{boundary}
\def\bsd/{\bdd/ \symdom/}
\def\bv/{boundary value}
\def\cf/{{\it cf}\.}
\def\Cf/{{\it Cf}\.}
\def\charr/{character}
\def\coeff/{coefficient}
\def\comm/{commutative}
\def\cpct/{compact}
\def\compl/{complex}
\def\comp/{complex}
\def\Comp/{Complex}
\def\conf/{conformal}
\def\conj/{conjugate}
\def\conn/{connect}
\def\cont/{continuous}
\def\conv/{converge} 
\def\convc/{convergence}
\def\convt/{convergent}
\def\convx/{convex}
\def\coord/{coordinate}
\def\lcoord/{local coordinate}
\def\Corr/{Corresponding}
\def\corr/{corresponding}
\def\corrd/{correspond}
\def\cov/{covariant}
\def\decomp/{decomposition}
\def\deco/{decompose}
\def\diff/{different} 
\def\Diff/{Different} 
\def\dimn/{dimension} 
\def\distr/{distribution} 
\def\div/{diverge} 
\def\dom/{domain}
\def\eg/{\hbox{\it e.g}\.}
\def\eigenf/{eigen\-\func/}
\def\eigensp/{eigen\-space}
\def\eigenv/{eigen\-value}
\def\eq/{equation}
\def\equa/{equation}
\def\de/{\diff/ial \equa/}
\def\do/{\diff/ial operator}
\def\ode/{ordinary \de/}
\def\pde/{partial \de/}
\def\pdo/{partial \diff/ial operator}
\def\psdo/{pseudo \diff/ial operator}
\def\fin/{finite}
\def\Ex/{\it Example.\ \normal}
\def\Exnr#1/{\it Example #1.\ \normal}
\def\foll/{follow}
\def\follg/{following}
\def\Follg/{Following}
\def\func/{function}
\def\Func/{Function}
\def\Fonc/{Fonc\-tion}
\def\fonc/{fonc\-tion}
\def\Funk/{Funk\-tion}
\def\funk/{Funk\-tion}
\def\gen/{general}
\def\har/{harmonic}
\def\Hint/{\it Hint. \normal}
\def\hist/{historic}
\def\histcl/{historical}
\def\hol/{holo\-morphic}
\def\homog/{ho\-mo\-ge\-ne\-ous}
\def\hyp/{hyper\-bolic}
\def\hyperg/{hyper\-geometric}
\def\ie/{\hbox{\it i.e.}}
\def\iff/{if and only if}
\def\ineq/{inequality}
\def\infra/{{\it inf\-ra}}
\def\ultra/{{\it ult\-ra}}
\def\Inpart/{In particular}
\def\inpart/{in particular}
\def\instof/{instead of}
\def\interps/{interpolation space}
\def\interp/{interpolation}
\def\Interp/{Interpolation}
\def\interpr/{Interpretation}
\def\Intr/{Introduction}
\def\intv/{interval}
\def\inv/{invariant}
\def\invc/{invariance}
\def\Iowords/{In other words}
\def\iowords/{in other words}
\def\ipr/{inner product}
\def\irred/{irreducible}
\def\lb/{line bundle}
\def\lin/{linear}
\def\lhs/{left hand side}
\def\rhs/{right hand side}
\def\loc/{local}
\def\math/{mathematic} 
\def\mathcn/{\math/ian}
\def\manif/{manifold}
\def\meas/{measure}
\def\measl/{measurable}
\def\mero/{mero\-morphic}
\def\mon/{monomial}
\def\monog/{monogenic}
\def\mult/{multiple}
\def\multy/{multiply}
\def\multn/{multiplication}
\def\nas/{necessary and sufficient}
\def\nbd/{neighborhood}
\def\neg/{negative}
\def\nondeg/{nondegenerate}
\def\Oohand/{On the other hand}
\def\oohand/{on the other hand}
\def\Oonhand/{On the one hand}
\def\oonhand/{on the one hand}
\def\oper/{operator}
\def\orth/{ortho\-gonal}
\def\orthon/{ortho\-normal}
\def\otoh/{on the other hand}
\def\quat/{quaternion}
\def\pp/{\hbox{a. e.}}
\def\psh/{plurisubharmonic}
\def\pol/{polynomial}
\def\pot/{potential}
\def\pos/{positive}
\def\princ/{principle}
\def\prob/{probability}
\def\proj/{projective}
\def\projn/{projection}
\def\Proof/{\it Proof:\normal}
\def\Rem/{\it Remark\normal}
\def\Remnr#1/{\it Remark\ \normal #1. }
\def\rep/{representation}
\def\meta/{metaplectic representation}
\def\repr/{reproducing}
\def\reprker/{reproducing kernel}
\def\resp/{respective} 
\def\resply/{respectively}
\def\restr/{restriction}
\def\sa/{self-adjoint}
\def\st/{such that}

\def\sol/{solution}
\def\ru/{space}
\def\sph/{spherical}
\def\ssp/{sub\ru/}
\def\sym/{symmetric}
\def\Sym/{Symmetric}
\def\symb/{symbol}
\def\symbc/{symbolic}
\def\symdom/{\sym/ domain}
\def\symp/{symplectic}
\def\Theor#1/{\fet Theorem #1.\ \normal}
\def\Lem#1/{\fet Lemma #1.\ \normal}
\def\Lemma/{\fet Lemma.\ \normal}
\def\topl/{topology}
\def\topll/{topological}
\def\transf/{transform}
\def\transl/{translation}
\def\transfn/{transformation}
\def\transv/{transvectant}
\def\trig/{trigonometric}
\def\tril/{trilinear}
\def\trilf/{trilinear form}
\def\uhp/{upper halfplane}
\def\uhs/{upper halfspace}
\def\vb/{vector bundle}
\def\vf/{vector field}
\def\vsp/{vector space}
\def\wrt/{with respect to}
\def\Wlog/{Without loss of generality}
\def\a{\alpha}

\def\sig{\sigma}

\def\Ab/{Abel}
\def\Ban/{Banach}
\def\Bansp/{\Ban/ space}
\def\Belt/{Bel\-tra\-mi}
\def\Berg/{Berg\-man}
\def\Bern/{Ber\-nou\-lli}
\def\Berz/{Berezin}
\def\Bess/{Bessel}
\def\Cart/{Car\-tan}
\def\Cay/{Cay\-ley}
\def\CG/{Clebsch-Gordan}
\def\Cl/{Clifford}
\def\CR/{Cauchy-Rie\-mann}
\def\Dir/{Dirichlet}
\def\Eucl/{Euclide}
\def\F/{Fourier}
\def\Hank/{Hankel}
\def\Hankf/{\Hank/ form}
\def\Herm/{Hermite}
\def\Hilb/{Hilbert}
\def\Hilbs/{Hilbert space}
\def\Hilbsp/{Hilbert space}
\def\HS/{Hilbert-Schmidt}
\def\Lag/{La\-grange}
\def\Lap/{La\-place}
\def\LapBelt/{\Lap/-\Belt/}
\def\Leb/{Lebesgue}
\def\Marc/{Mar\-cin\-kie\-wicz}
\def\Moeb/{Moebius}
\def\Moebt/{Moebius transformation}
\def\Moebtransfn/{Moebius transformation}
\def\Pla/{Plan\-che\-rel}
\def\Poin/{Poin\-car\'e}
\def\Riem/{Rie\-mann}
\def\Riemn/{\Riem/ian}
\def\psRiemn/{pseudo-\Riem/ian}
\def\Riems/{Rie\-mann surface}
\def\Schroe/{Schr\"odinger}
\def\Weier/{Weier\-strass}
\def\anal/{analytic}
\def\bsd/{bounded symmetric domain  }
\def\bdd/{bounded}
\def\calc/{calculation}\def\conj{conjugate}
\def\calci/{calculating}\def\eg{e.g.}
\def\conj/{conjugate}
\def\deco/{decomposition}
\def\eg/{e.g.}
\def\fct/{function}
\def\gp/{group}
\def\hw/{highest weight}
\def\hwv/{highest weight vector}
\def\hwvs/{highest weight vectors}
\def\lw/{lowest weight}
\def\lwv/{lowest weight vector}
\def\lwvs/{lowest weight vectors}
\def\hds/{holomorphic discrete series}
\def\iff/{if and only if}
\def\inv/{invariant}
\def\irrde/{irreducible decomposition}
\def\meas/{measure}
\def\transf/{transform}
\def\rep/{representation}
\def\resp/{respectively}
\def\inters/{intertwines}
\def\interg/{intertwining}
\def\meta/{metaplectic representation}
\def\qu/{quaternion}
\def\rep/{representation}
\def\symdom/{ symmetric domain}
\def\st/{such that}
\def\shd/{subhead}
\def\transf/{transform}
\def\wrt/{with respect to}

\def\Norm#1#2#3{\Vert#1\Vert^{#3}_{{#2}¥}}

\begin{document}
\title[Laplacians on Cauchy-Riemann complexes and and harmonic forms
]
{Laplacians on quotients of Cauchy-Riemann complexes and 
Szeg\"o map for $L^2$-harmonic forms 
}
\author{Bent \O{}rsted and Genkai Zhang}
\address{B\O{}: Department of Mathematics, Aarhus University,
Denmark
}\address{GZ:Department of Mathematics, Chalmers University of Technology and
Gothenburg University
 , Gothenburg, Sweden}
\email{orsted@imada.sdu.dk, genkai@math.chalmers.se}

\dedicatory{Dedicated to Professor Gerrit van Dijk on
the occasion of his 65th birthday}

\thanks{Research by Genkai Zhang supported by the Swedish
Science Council (VR) }

\def\GG{\mathbb G}
\def\KK{\mathbb K}
\begin{abstract}

We compute the spectra of the Tanaka type Laplacians
$\square=\bar \partial_{\mathcal Q}^\ast\bar \partial_{\mathcal Q}
+\bar \partial_{\mathcal Q} \bar \partial_{\mathcal Q}^\ast $ and 
$\triangle =\partial_{\mathcal Q}^\ast \partial_{\mathcal Q}
+ \partial_{\mathcal Q} \partial_{\mathcal Q}^\ast $
on the Rumin complex $\mathcal Q$, a quotient
of the tangential Cauchy-Riemann complex on the unit
sphere $S^{2n-1}$ in $\mathbb C^n$.
We prove that Szegö map
is a unitary operator from
a subspace of $(p, q-1)$-forms on the sphere 
defined by the operators $\triangle$
and the normal vector field
onto the space of
 $L^2$-harmonic
$(p, q)$-forms on the unit ball.
Our results  generalize earlier result of Folland.
\end{abstract}

\maketitle








\def\mP{\mathcal P}
\def\mF{\mathcal F}
\def\SdN{\mathcal S^N}
\def\eSdN{\mathcal S_N^\prime}
\def\Gc{G_{\mathbb C}}
\def\rqs{G\backslash K_{\mathbb C}}

\baselineskip 1.25pc

\section{Introduction}

The  Bergman space of holomorphic functions
on the unit ball $B$ of $\mathbb C^n$ plays
an important role both in representation theory
of the group $SU(n, 1)$
and complex analysis. The holomorphic functions
in the space can be identified with
square integrable  harmonic $(n, 0)$-forms
on $B$ equipped with the Bergman metric.
One may naturally consider the $(p, q)$-forms
on the unit ball.
It has been proved by Donnelly and Fefferman
\cite{Don-Fef} that the space of $L^2$ harmonic
$(p, q)$ forms on a strongly pseudo-convex domain
is non-empty precisely when $p+q=n$; 
the case of unit ball  was proved earlier
in representation theory of semisimple Lie groups
and is of much importance 
for the general strongly pseudo-convex domains.
The study of
$L^2$ harmonic $(p, q)$-forms on the unit ball 
is already quite involved, partly because the existence
theorem does not give explicit construction. In the case
of the real ball in $\mathbb R^{2n}$ with the hyperbolic metric
the $L^2$-harmonic $n$-forms (with respect to
the corresponding Laplacian) are Euclidean
harmonic forms. One can thus consider the
boundary value map of the harmonic forms on the unit sphere,
its dual map gives then the Szegö map from $(n-1)$-form
to harmonic $n$-forms on the real ball. 
The Szegö map is special case of the general
Poisson transform mapping the forms on the sphere
to eigenforms on the ball of invariant differential operators.
In 
\cite{Gaillard-cmh-86} 
Gaillard computed explicitly the Poisson
transform from $(p-1)$-forms on the boundary
to harmonic (not necessarily $L^2$) $p$-forms on the ball.
Branson 
\cite{Branson-korea}, Chen \cite{Chen-Kthy}, and Lott \cite{Lott-cmh-cmh-00}
have found  formulas expressing the $L^2$-norm
of the harmonic forms in terms of their boundary values
and in terms of the Laplacian operators on the boundary.
These formulas turn out to be quite important
both in conformal geometry, $K$-theory and in representation
theory, see loc. cit. 
Part of our purpose is to find a corresponding formula
for the unit complex ball.

In a remarkable paper 
\cite{Julg-Kasparov} 
Julg and Kasparov 
give a geometric construction of $L^2$-harmonic $(p, q)$-forms,
$p+q=n$,
on the unit ball and of the Szeg\"o map from forms
on the sphere into harmonic forms. On the sphere $S$ 
there is the Cauchy-Riemann complex 
defined as space of $(p, q-1)$-forms
modulo the ideal generated by
the contact form $\bar \tau=\sum_{j=1}^n \bar z_j d
z_j$.
However it is proved in 
\cite{Julg-Kasparov} that, roughly speaking,
 the boundary forms
of $L^2$ harmonic $(p, q)$-forms are
in the so-called Rumin complex $\mathcal Q$, which is
the quotient of Cauchy-Riemann complex
modulo the forms $\tau$ and $d\tau
=\sum_{j=1}^n  d\bar z_j\wedge dz_j$. In the middle
of the Rumin complex is the space of $n-1$-forms
and there defined a second-order
Rumin operator in place of the
the first order boundary C-R operator $\bar\partial_b$,
the boundary value of harmonic forms
are in certain Sobolev-type space defined by the
Rumin Laplacian which is of $4$-order.
In the present paper
we will use their results
to prove a formula expressing the $L^2$-norm of
the harmonic forms in terms
certain Sobolev norm of their boundary values
defined in terms of the Tanaka-type Laplacians
$\square_{\mathcal Q} =\bar \partial_{\mathcal Q}^\ast 
\bar \partial_{\mathcal Q} +
\bar \partial_{\mathcal Q} \bar \partial_{\mathcal Q}^\ast $
and $\triangle_{\mathcal Q} = \partial_{\mathcal Q}^\ast 
 \partial_{\mathcal Q} +
 \partial_{\mathcal Q} \partial_{\mathcal Q}^\ast $
on the complex $\mathcal Q$
on the unit sphere. For that purpose we compute
the spectrum of the  Laplacians on the
Rumin complex of $(p, q)$-forms
for general  $(p, q)$. 
For $(0, q)$-forms
the operator $\square_{\mathcal Q}$ coincides with
the usual C-R Laplacian $\square_{b}$
and its 
 spectrum has been found earlier
by Folland \cite{Folland-tams-72}.
 So our results generalize
that of Folland.

The $L^2$-space of differential
forms in the Rumin complex on the sphere
is a natural
model for the induced representations
of the group $SU(1, n)$ from certain
finite-dimensional representations
of the maximal parabolic subgroup.
There arise many other geometrically interesting differential
operators
and unitary representations such as complementary series,
see e.g.\cite{Branson-tams95}.
We hope that our results here will shed light
on studying those subjects.

The paper is organized as follows. In Section 2 we 
present the context of our problems and recall
the known results. In Section 3 we compute
the the eigenvalues of the Laplacians
$\square_b$ and $\triangle_b$
on 
certain subspaces of $(p, q-1)$ forms.
The unitary structure of the space of $L^2$
harmonic $(p,q)$-forms is given in Section 4. 

We would like to thank the organizers
of the  conference ``Harmonic Analysis and Homogeneous
Spaces'', Lorentz Center, Leiden University, August 23-27, 2004,
for the invitation.
Thanks are also  due to Thomas Branson for some
helpful discussions and to Miraslav Engli\v{s} for reading
an earlier version of this paper.

\section{Cauchy-Riemann complexes and their quotients}

For some notational convenience  we
will consider $(p, q-1)$-forms on the unit sphere $S$.
We assume
throughout the text that $n>1$, $1\le p, q\le n$.

We recall first briefly some well-known
facts about Cauchy-Riemann complex
of differential forms on $S$.

Let $\Ome^{p, q-1}$
be the space of smooth
differential $(p, q-1)$-forms on some neighborhood 
of the unit sphere $S$. Recall by definition
that $(p, q)$-form
is a section of the vector bundle
$\wedge^{q} {T^{(0, 1)}}^\prime
\otimes \wedge^{p} {T^{(1, 0)}}^\prime$
over some neighborhood of $S$. For  simplicity we denote also
$\Ome^{p, q-1}_z= \wedge^{q-1} {T_z^{(0, 1)}}^\prime
\otimes \wedge^{p} {T_z^{(1, 0)}}^\prime$
the fiber space at each $z$. We fix also the usual convention
that the $(1, 0)$-covectors in ${T_z^{(1, 0)}}^\prime$
and the $(0, 1)$-covectors
in $ {T_z^{(0, 1)}}^\prime$ are anticommuting.

 We let 
$$
\tau=\sum_{j=1}^n \bar{z_j}dz_j, \quad
\bar \tau=\sum_{j=1}^n {z_j}d\bar z_j, \quad
$$
be the  $(1, 0)$- and $(0, 1)$-forms defining
the contact structure on $S=\partial B$.
 We let
 $\Omega^{p, q-1}/(\bar \tau)$
be the quotient
by
the ideal $ (\bar \tau)$ generated by $\bar\tau$. Since
$\bar\partial\bar\tau=0$, there is
an induced operator, denoted by $\bar\partial_b$
on the quotient
 $\Omega^{p, q-1}/(\bar \tau)$ and it forms the C-R
complex.

For our purpose we need to consider some further quotients.
Let
 $\mathcal Q^{p, q-1}= \Ome^{p, q-1}/(\tau, \bar \tau, d\tau)$ be the quotient
of $\Ome^{p, q-1}$ by the ideal generated by
$\tau$,  $\bar \tau$ and $d\tau$; the ideal is clearly
invariant under $\partial$ and $\bar\partial$,  the induced
operators  on the quotient
$\mathcal Q^{p, q-1}$
will  be denoted by $\partial_{\mathcal Q}$ and 
$\bar\partial_{\mathcal Q}$. The spaces
$\sum_{p+q-1=j}\mathcal Q^{p, q-1}$ for $j\le n-1$ 
and the induced operator of $d=\bar\partial +\partial$
can be prolonged to a complex, also called Rumin complex,
so that the cohomology so obtained is
the usual de Rham cohomology, see 
\cite{Julg-Kasparov} and references therein.

Consider the unitary
group $U(n)$ acting natural on tangent space $\mathbb C^n
=\sum_{j=1}^n\mathbb C\frac{\partial}{\partial z_j}$
of $B$ 
by definition, and its dual action
on  the differential forms
$\sum_{j=1}^n \mathbb C dz_j$ of degree $1$
with constant coefficients. 
We fix the unique $U(n)$-invariant
point-wise inner product $\langle \cdot, \cdot\rangle$
on the differential forms
at each point normalized
so that
$
 d\bar z_{j_1}\wedge \cdots \wedge  d\bar z_{j_q}
\wedge
 d z_{i_1}\wedge \cdots \wedge d z_{i_p}$,
$1\le i_1 <
\cdots < i_p\le n$,
$1\le j_1 <\cdots < j_q$
form an orthonormal basis of different $(p,q)$ forms 
at each $z\in S$.
We
define the $U(n)$-invariant
inner product on 
 $\Ome^{p, q-1}$ by
$$
(\alpha, \beta)=\int_S \langle \alpha(z), \beta(z)\rangle 
d\sig(z)
$$
where $d\sig(z)$ is the normalized area measure
on $S$.

We need further some important
operators on differential forms. For $\a\in \Ome^{p, q-1}$   we denote
$L_{\alpha}$
the corresponding multiplication operators by $\a$,
$L_{\alpha}\beta =\alpha \wedge \beta$. For a holomorphic
tangent or anti-holomorphic tangent vector $X$ in $\mathbb C^n$
we let $i_X$ be the standard evaluation operator
see e.g. \cite[Chapter 0, Section 7]{Griffiths-Harris}. In particular
the operators $\partial, L_{dz_j}, i_{\partial_j}$
anti-commutes with the operators
$\bar\partial, L_{d\bar z_j}, i_{\bar\partial_j}$.
Let 
$$
T=\sum_{j=1}^n z_j\partial_j, \,\, 
\bar T=\sum_{j=1}^n \bar z_j\bar\partial_j, \,\, 
$$
and
$$
E=\sum_{j=1}^n \partial_j\otimes \bar \partial_j
$$
and $i_E$ the corresponding evaluation operator
$$
i_E=\sum_{j=1}^n i_{\partial_j} i_{\bar \partial_j}.
$$
Viewed as operators on $\Ome^{p, q-1}$ with the 
point-wise inner product $\langle \alpha(z), \beta(z)\rangle $
for each  point $z$
their adjoint operators are 
\begin{equation}
  \label{eq:duality}
(L_{d\tau})^\ast=i_E,\quad (L_{\tau})^\ast=i_T, \quad
(L_{\bar\tau})^\ast=i_{\bar T}.  
\end{equation}
This is then also true with respect to the
inner product $(\alpha, \beta)$.

The space  $\mathcal Q^{p, q-1}=\Ome^{p, q-1}/{(\tau, \bar \tau, d\tau)}$ can  be realized
as the subbundle $\Ker i_{\bar T}\cap \Ker i_{T}\cap \Ker i_E$ on
$\Ome^{p, q-1}$ and we will fix this realization
throughout our paper. The induced
C-R operators $\partial_{\mathcal Q}
$ and
$\bar \partial_{\mathcal Q}$ will be, under
this identification, the differentiations
$\partial$ and $\bar \partial$,
followed by the projection onto
the sub-bundle.
We put
\begin{equation}
  \label{eq:zeta-xi}
\zeta_j:=\bar\partial_b \bar z_j =
 \bar\partial_{\mathcal Q} \bar z_j=  d\bar z_j -\bar z_j \bar \tau, \quad
\xi_j=\partial_{\mathcal Q} z_j  =dz_j - z_j  \tau;  
\end{equation}
 see e.g. Folland \cite{Folland-tams-72}. 

The next lemma is elementary.

\begin{lemm+} The differential
form $d\tau-\bar \tau\wedge \tau$
is $T$ and $\bar T$-primitive, namely
it is annihilated by $i_{T}$
and $i_{\bar T}$. Moreover the forms 
$\zeta_j$ and 
$\xi_j$
can also be obtained by
$$
\zeta_j=-i_{\partial_j}(d\tau -\bar\tau \wedge \tau), \qquad
\xi_j=i_{\bar \partial_j}(d\tau -\bar\tau \wedge \tau).
$$
\end{lemm+}

We need also some Hodge type commutation relations;
the first equality below is  well known 
(see e.g. \cite[Chapter 0, Section 7]{Griffiths-Harris}), the other
can be proved by routine computation.

\begin{lemm+}
Assume $p+q\le n$. The following commutation
relations hold,  viewed as operators
acting on the space $\Omega^{p, q-1}$, 
$$[i_E, L_{d\tau}]=p+q
$$
$$[i_T, L_{d\tau}]=L_{\bar \tau}, \quad
[i_{\bar T}, L_{d\tau}]=L_{\tau}
$$
and
$$
i_T L_{\tau} + L_{\tau}i_T  = |z|^2, \quad
i_{\bar T} L_{\bar \tau} + L_{\bar \tau}i_{\bar T}  = |z|^2.
$$
\end{lemm+}

By taking adjoint we get then
\begin{equation}
  \label{comm-ie-lt}
[i_E, L_{\tau}]=-i_{T}, \quad
[i_{E}, L_{\bar \tau}]=i_{\bar T}.
\end{equation}
It follows further  that
\begin{equation}
  \label{comm-1}
[i_T, L_{d\tau-\bar \tau \wedge \tau}]=0,
\end{equation}
and
\begin{equation}
  \label{eq:comm}
[i_E-i_{\bar T} i_{T}, L_{d\tau-\bar \tau\wedge \tau}]
=(n-1)-\sum_{j=1}^n L_{\zeta_j}i_{\partial_j}
-\sum_{j=1}^n L_{\xi_j}i_{\bar \partial_j}.
\end{equation}

The following decomposition, which
will not be used in this paper,
follows
easily from our previous commutation relation, and
it might be of some independent interests; it
can be proved easily and we omit the proof.

\begin{coro+}Realizing
the quotient
$\mathcal Q^{p, q-1}$ 
as
the space 
$\Ker i_{T}\cap \Ker i_{\bar T}\cap \Ker {i_E}$
on $\Ome^{p, q-1}$,
we have
$$
\Ome^{p, q-1}=\mathcal Q^{p, q-1} \oplus
(L_{\tau}\Ome^{p-1, q-1} + L_{\bar \tau} \Ome^{p-1, q-1})
 \oplus L_{d\tau-\bar \tau\wedge \tau}\Ome^{p-1, q-2}
$$
\end{coro+}

To state the next result we fix the space
$\mathfrak{h}=\mathbb C H_1+\cdots +\mathbb C H_n$,  of diagonal
matrices as a Cartan subalgebra of the Lie algebra
$\mathfrak{u}(n)^{\mathbb C}= \mathfrak{gl}(n, \mathbb C)$, with
$H_j=\text{diag}(0, \cdots, 0, 1, 0, \cdots, 0)$, $j=1, \cdots, n$.
Let $\{\varepsi_j\}$ be the basis of $\mathfrak{h}^\prime$ dual to
$\{H_j\}$, and we fix an ordering of roots of 
$\mathfrak{h}$ in $\mathfrak{gl}(n, \mathbb C)$ so that
$\varepsi_1 > \cdots >\varepsi_n$.

\begin{prop+}\label{m-free}
Assume $p+q\le n$. As a representation
space of $K=U(n)$ the $L^2$-space
of  $(p, q-1)$-forms in the quotient $\mathcal Q^{p, q-1}=\Ome^{p, q}/(\tau, \bar \tau, d\tau)$
is decomposed 
with multiplicity free as follows
$$
L^2(\mathcal Q^{p, q-1}
)
=\sum_{m, l\ge 0}\!\!\!\!\!{\phantom{|}^{\phantom{|}}}^\oplus(\Phi^{m, l}_{q-1, p}\oplus
\Phi^{m, l}_{q-1, p+1}
\oplus
\Phi^{m, l}_{q, p}
\oplus
\Phi^{m, l}_{q, p+1}),
$$
where $\Phi^{m, l}_{r, s}$
 are 
of  highest weight 
$$(m, 0, \cdots, 0, -l)+(\overbrace{1, \cdots 1}^r, 0, \cdots, 0, 
\overbrace{-1,
\cdots, -1}^s)$$ 
with  $r$ many $1$'s and $s$ many $-1$'s,
for $(r, s)=(q-1, p), (q-1, p+1), (q, p), (q, p+1)$;
the last summand will not appear when $p+q=n$.
\end{prop+}

\begin{proof} This follows from the Frobenius reciprocity
and is stated in \cite{Julg-Kasparov}  (without proof).
We give a brief proof here. View
the space $S$ as a homogeneous space $S=U(n)/U(n-1)$
where $U(n-1)$ is the subgroup
keeping the base point
 $z_0=(1, 0, \cdots, 0)\in S$ fixed.
The bundle $\mathcal Q^{p, q-1}$ is a $U(n)$-homogeneous bundle,
and the fiber space at
the base point is spanned by the vectors
$ d\bar z_{i_1}\wedge \cdots \wedge  d\bar z_{i_{q-2}}
\wedge
 d z_{j_1}\wedge \cdots \wedge d z_{j_{p-1}}$,
$2 \le i_1 <
\cdots < i_{q-1}\le n$,
$2\le  j_1 <\cdots < j_{p-1}\le n$ modulo
the ideal generated by the form
$d\bar z_2\wedge d z_2+\cdots+
d\bar z_{n}\wedge d z_n$
since $\tau(z_0)=dz_1$ and $\tau(z_0)=d\bar z_1$,
and it forms an irreducible representation $\omega_0$
of $U(n-1)$, by e.g. the Lefschetz decomposition 
\cite{Griffiths-Harris} and by
the tensor product decomposition 
of representations of $U(n-1)$ \cite{Ze}.
 The multiplicity
of a representation $\mu$
of $U(n)$ appearing
in 
$L^2(\mathcal Q^{p, q-1}
)$
is the same as
the multiplicity of $\omega_0$
in $\mu$. The rest follows
then by the known
results on branching under $U(n-1)$
of representations of $U(n)$,
see e.g. \cite{Ze}.
\end{proof}

\section{Spectrum of the  Laplacians
$\square_{\mathcal Q} =\bar \partial_{\mathcal Q}^\ast 
\bar \partial_{\mathcal Q}
+\bar \partial_{\mathcal Q}\bar \partial_{\mathcal Q}^\ast 
$
and $\triangle_{\mathcal Q} = 
\partial_{\mathcal Q}^\ast 
 \partial_{\mathcal Q}
+\partial_{\mathcal Q}\partial_{\mathcal Q}^\ast $}

In this section we compute the spectrum of the Tanaka-type
Laplacians $\square_{\mathcal Q}$
and $\triangle_{\mathcal Q} $ in the space $L^2({\mathcal Q}^{p, q-1})
$, by computing their eigenvalues on each irreducible
subspace.
To simplify the notation we will suppress the upper
indices $(m, l)$.

\begin{lemm+}The  highest weight vectors 
$\alpha, \beta, \gamma, \delta$
in the 
spaces
$\Phi^{m, l}_{q-1, p}, \Phi^{m, l}_{q-1, p+1}, \Phi^{m, l}_{q, p}$
and  $\Phi^{m, l}_{q, p+1}$ 
 are given as follows:
$$
 \a_{q-1, p}= \bar z_1^{m} z_n^{l} \a_{q-1, p}^\prime 
$$
with
\begin{equation*}
\begin{split}
\a_{q-1, p}^\prime =&\zeta_1\wedge \cdots \wedge \zeta_{q-1}
\wedge \xi_{n-p+1}\wedge \cdots \wedge \xi_{n}\\
&\quad -\frac 1{n-1}(d\tau-\bar \tau \wedge \tau)
\wedge i_{E}(\zeta_1\wedge \cdots \wedge \zeta_{q-1}
\wedge \xi_{n-p+1}\wedge \cdots \wedge \xi_{n}),
\end{split}
\end{equation*}
$$
\beta_{q-1, p+1}=\bar z_1^{m} z_n^{l}
\zeta_1\wedge \cdots \wedge \zeta_{q-1}
\wedge \left(
\sum_{i=n-p}^n (-1)^{i-(n-p)}z_i\xi_{n-p}\wedge
 \cdots \widehat{\xi_i} \cdots
\wedge \xi_{n}\right),
$$
$$
\gamma_{q, p}=\bar z_1^{m} z_n^{l}
\left(\sum_{j=1}^q(-1)^{j-1}\bar z_j
\zeta_1\wedge \cdots \widehat{\zeta_j}\cdots\wedge \zeta_{q}
\right)\wedge \xi_{n-p+1}\wedge \cdots\wedge \xi_{n}
$$
and respectively
\begin{equation*}
\begin{split}
&\qquad \delta_{q, p+1}\\
&=\bar z_1^{m} z_n^{l}
\left(\sum_{j=1}^q(-1)^{j-1}\bar z_j
\zeta_1\wedge \cdots \widehat{\zeta_j}\cdots \wedge\zeta_{q}
\right)\wedge \left(
 \sum_{i=n-p}^n (-1)^{i-(n-p)}z_i\xi_{n-p}
\wedge \cdots \widehat{\xi_i} \cdots
\wedge \xi_{n}\right).
\end{split}
\end{equation*}
\end{lemm+}

\begin{proof} It is easy to prove that
the given forms are of the respective weights. By the multiplicity
free result in 
Proposition \ref{m-free} we need only to prove that
they are in $\mathcal Q$.
We consider first the forms $\beta, \gamma, \delta$. 
It follows from \cite[Theorem 5]{Folland-tams-72}
that
\begin{equation}
  \label{eq:folland-1}
\omega_1(q):=\sum_{j=1}^q(-1)^{j-1}\bar z_j
\zeta_1\wedge \cdots \widehat{\zeta_j}\cdots \wedge\zeta_{q}
=\sum_{j=1}^q(-1)^{j-1}\bar z_j
d\bar z_1\wedge \cdots \widehat{d\bar z_j}\cdots
\wedge d\bar z_{q},  
\end{equation}
and (by rewriting the formula there)
\begin{equation}
  \label{eq:folland-2}
\zeta_1\wedge \cdots \wedge \zeta_{q-1}
=d\bar z_1\wedge \cdots \wedge d\bar z_{q-1}
-\bar \tau\wedge \omega_1(q-1) 
\end{equation}
and they are highest weight vectors
in the C-R complex
 $\mathcal Q^{0, q}$ of $(0, q)$-
forms. Similarly
 we can prove that
\begin{equation}
  \label{eq:folland-3}
\begin{split}
 \omega_2(p):&=\sum_{i=n-p}^n (-1)^{i-(n-p)}z_i\, \xi_{n-p}
\wedge \cdots\wedge \widehat{\xi_i} \wedge\cdots
\wedge \xi_{n}\\
&=\sum_{i=n-p}^n (-1)^{i-(n-p)}z_i dz_{n-p}
\wedge \cdots \widehat{dz_i} \cdots
\wedge dz_{n}
\end{split}
\end{equation}
and
\begin{equation}
  \label{eq:folland-4}
\xi_{n-p+1}\wedge \cdots \wedge \xi_{n}
= dz_{n-p+1}\wedge \cdots \wedge 
dz_{n} - \tau\wedge \omega_2({p-1})
\end{equation}
are highest weight vectors
in the C-R complex
 $\mathcal Q^{p, 0}$ of $(p, 0)$-
forms. The forms
$\beta$,  $\gamma$ and $\delta$ 
are thus all annihilated by
$i_{T}$ and $i_{\bar T}$; moreover
we have $i_E(\delta)=0$ by a simple computation.
We  prove now that $i_E(\gamma)=0$. We have,
by the above formulas,
$$
\gamma=\omega_1(q)\wedge dz_{n-p+1}\wedge \cdots
\wedge dz_{n}-(-1)^{q-1}\tau\wedge \omega_1(q)\wedge \omega_2(p-1).
$$
The first term is clearly annihilated by $i_E$,
and so is also $\omega_1(q)\wedge \omega_2(p-1)$. 
Thus by the commutation relation 
  (\ref{comm-ie-lt}), we have
that up to a factor of $\pm 1$, $
i_E(\gamma)$ is
$$i_{\bar T}(\omega_1(q)\wedge \omega_2(p-1))=0$$
since each term $ \omega_1(q)$ and $ \omega_2(p-1)$
is annihilated by $i_{\bar T}$.
The proof of $i_E(\beta)=0$ is the same. 

Consider finally the form $\alpha$. 
We have that $i_T(\alpha)=0$,
$i_{\bar T}(\alpha)=0$, since
the forms $\zeta$'s and $\xi$'s
and $d\tau-\bar \tau\wedge\tau$ are
annihilated by $i_T$ and $i_{\bar T}$.
We prove now $i_E(\alpha)=0$.
Write 
$$
\alpha^\prime=\kappa -\frac 1{n-1}(d\tau-\bar \tau \wedge \tau)
\wedge i_E(\kappa)$$
with 
$$\kappa =\zeta_1\wedge\cdots \wedge\zeta_{q-1}\wedge
\xi_{n-p+1}\wedge\cdots \wedge\xi_{n}.
$$
Now 
\begin{equation*}
\begin{split}
i_E(\alpha^\prime)&
=i_E(\kappa) -\frac 1{n-1}i_E
\left((d\tau-\bar \tau \wedge \tau )
i_E(\kappa)\right)\\
&=i_E(\kappa) -
\frac 1{n-1}\left(i_E(d\tau-\bar \tau \wedge \tau)i_E(\kappa) +(d\tau-\bar \tau \wedge \tau)\wedge i_Ei_E(\kappa)\right)\\
&=i_E(\kappa) -
\frac 1{n-1}(n-1)i_E(\kappa) -\frac 1{n-1}
(d\tau-\bar \tau \wedge \tau)\wedge i_Ei_E(\kappa)\\
&=-\frac 1{n-1}
(d\tau-\bar \tau \wedge \tau)\wedge i_E
i_E(\kappa)
\end{split}
\end{equation*}
since $(i_E(d\tau)-i_E(\bar \tau \wedge \tau))=n-1$.  We claim
that
\begin{equation}
  \label{eq:ie2}
i_E
i_E(\kappa)=0.
\end{equation}
Indeed, 
the form  $\kappa$  is 
\begin{equation*}
\begin{split}
\kappa&=d\bar z_1\wedge\cdots \wedge d\bar z_{q-1}
\wedge
dz_{n-p+1}\wedge\cdots \wedge dz_{n}
-\bar \tau \wedge \omega_1\wedge dz_{n-p+1}\wedge\cdots \wedge dz_{n}
\\
&\qquad 
- (-1)^{q-1}\tau \wedge
d\bar z_1\wedge\cdots ...\wedge d\bar z_{q-1}\wedge \omega_2
+(-1)^ {q-1}   \bar \tau \wedge \tau\wedge \omega_1 
\wedge \omega_2 
\end{split}
\end{equation*}
with $\omega_1=\omega_1(q-1)$ and $\omega_2=\omega_2(p-1)$
defined in (3.3)-(3.4) ($q$ and $p$ replaced by $q-1$ and $p-1$ respectively).
From this we have 
\begin{equation*}
\begin{split}
i_E(\kappa)&-i_E(\bar \tau\wedge \omega_1\wedge
d z_{n-p+1}\wedge\cdots  d z_{n})
-(-1)^{q-1}i_E(\tau\wedge 
{d\bar z_1}
\wedge\cdots  d\bar z_{q-1}\wedge \omega_2)
\\
&\qquad + (-1)^{q-1}
i_E(\bar \tau\wedge \tau \wedge \omega_1
\wedge \omega_2).
\end{split}
\end{equation*}
Using the commutation relation 
  (\ref{comm-ie-lt})
and using the fact that
all the forms $d\bar z_1
\wedge\cdots  d\bar z_{q-1}\wedge \omega_2$, $\omega_1\wedge
d z_{n-p+1}\wedge\cdots d z_{n}$ 
and $\omega_1$ and $\omega_2$ are forms in 
$\mathcal Q$ we get 
\begin{equation}
\label{iEkappa}
i_E(\kappa)= (-1)^{q-1}\omega_1\wedge\omega_2,
\end{equation}
and (\ref{eq:ie2}) follows,
again by the formulas (3.1) and (3.3).

\end{proof}

The following lemma 
gives
the $L^2$-norms of the highest
weight vectors $\alpha$,
$\beta$, $\gamma$, and $\delta$. For
computational reasons we assume
$p+q<n$ for certain cases. If $p+q=n$
one can find the norm
by the same computation which we omit here.

\begin{lemm+}The following
formulas hold
\begin{equation*}
\begin{split}
& \Vert\alpha_{q-1, p}\Vert^2 
=\frac{m!l!}{(n)_{m+l+2}}\\
&\qquad \times \left((l+n-q+1)(m+n-p)
-\frac1{n-1}(l+p)(m+q-1)\right),  \quad (p+q <n)  
\end{split}
\end{equation*}
$$
\Vert\beta_{q-1, p+1}\Vert^2=\frac{m!l!}{(n)_{m+l+2}}(l+p+1)(l+n-q),
$$
$$
\Vert\gamma_{q, p}\Vert^2=\frac{m!l!}{(n)_{m+l+2}}(m+q)(m+n-p+1),
$$
$$
\Vert\delta_{q, p+1}\Vert^2=\frac{m!l!}{(n)_{m+l+2}}(l+p+1)(m+q).
\quad (p+q<n)
$$
\end{lemm+}
\begin{proof}
We consider first the forms  $\beta$, $\gamma$ and $\delta$.
The pointwise inner products $\langle \beta, \beta\rangle$,
$\langle \gamma, \gamma\rangle$, and 
$\langle \delta, \delta\rangle$
can be computed using 
  the equalities (\ref{eq:folland-1})-(\ref{eq:folland-4})
and the orthogonality
of the forms $d\bar z_{j_1}\wedge \cdots
d\bar z_{j_{q-1}}\wedge d z_{i_1}\wedge \cdots
d\bar z_{i_{p}}$ (see also  \cite{Folland-tams-72}, though
a different normalization was used), we have
$$
\langle \beta, \beta\rangle
=|\bar z_1^m z_n^l|^2(\sum_{i=n-p}^n |z_i|^2)
(\sum_{j=q}^{n} |z_j|^2),
$$
$$
\langle \gamma, \gamma\rangle
=|\bar z_1^m z_n^l|^2 
(\sum_{i=1}^{n-p} |z_i|^2)
(\sum_{j=1}^{q} |z_j|^2),
$$
$$
\langle \delta, \delta\rangle
=|\bar z_1^m z_n^l|^2 (\sum_{j=1}^{q} |z_j|^2)
(\sum_{i=n-p}^{n} |z_i|^2).
$$
Their $L^2$-norms can be computed then by integration
over $S$, and we omit the details.

We compute now the norm of $\alpha$. Adapting
the notation in the proof of Lemma 3.1, we have
\begin{equation*}
\begin{split}
\langle \alpha^\prime, \alpha^\prime
\rangle
&=\langle \kappa, \kappa
\rangle
-\frac 1{n-1}\langle \kappa, (d\tau-\bar \tau\wedge\tau)\wedge i_E\kappa
\rangle
-\frac 1{n-1}\langle (d\tau-\bar \tau\wedge\tau)\wedge i_E\kappa,
\kappa
\rangle \\
&\qquad +
\frac 1{(n-1)^2}\langle (d\tau-\bar \tau\wedge\tau)\wedge i_E\kappa,
(d\tau-\bar \tau\wedge\tau)\wedge i_E\kappa\rangle.
\end{split}
\end{equation*}
The inner product in the second term is
$$
\langle \kappa, (d\tau-\bar \tau\wedge\tau)\wedge i_E\kappa
\rangle
=\langle (i_E-i_Ti_{\bar T})\kappa, i_E\kappa
\rangle
=\langle i_E\kappa, i_E\kappa
\rangle
$$
since $\kappa $ is in $\Ker i_T$. So 
the third term is also $\langle i_E\kappa, i_E\kappa
\rangle$. The inner product
in the last term is
\begin{equation*}
\begin{split}
&\qquad \langle (d\tau-\bar \tau\wedge\tau)\wedge i_E\kappa,
(d\tau-\bar \tau\wedge\tau)\wedge i_E\kappa
\rangle
=\langle(i_E-i_Ti_{\bar T}) \left((d\tau-\bar \tau\wedge\tau)\wedge i_E\kappa\right),
i_E\kappa
\rangle \\
&=\langle i_E\left((d\tau-\bar \tau\wedge\tau)\wedge i_E\kappa\right),
i_E\kappa
\rangle
-\langle i_T\left(i_{\bar T}((d\tau-\bar \tau\wedge\tau)\wedge i_E\kappa)\right),
i_E\kappa
\rangle.
\end{split}
\end{equation*}
But  by (\ref{comm-1})
$$
i_{\bar T}((d\tau-\bar \tau\wedge\tau)\wedge i_E\kappa)
=(d\tau-\bar \tau\wedge\tau)\wedge i_{\bar T}i_E\kappa
$$
and $i_{\bar T}i_E\kappa=i_Ei_{\bar T}\kappa=0$. So 
this inner product is
\begin{equation*}
\begin{split}
  \langle i_E((d\tau-\bar \tau\wedge\tau)\wedge i_E\kappa),
i_E\kappa
\rangle &=(n-1)\langle i_E\kappa,
i_E\kappa
\rangle
+\langle (d\tau-\bar \tau\wedge\tau)\wedge(i_E
i_E\kappa),
i_E\kappa
\rangle \\
=(n-1)\langle i_E\kappa,
i_E\kappa
\rangle
\end{split}
\end{equation*}
by  (\ref{eq:ie2}).
We have then
\begin{equation}
\label{inn-a}
\langle \alpha, \alpha
\rangle  =|\bar z_1^m  z_n^l |^2
\langle \kappa, \kappa
\rangle  -\frac 1{n-1}
|\bar z_1^m  z_n^l |^2
\langle i_E\kappa,
i_E\kappa
\rangle.
\end{equation}
The point-wise norm squared of $\kappa$ is
$$
\langle \kappa, \kappa\rangle
=(\sum_{i=1}^{n-p}|z_i|^2)
(\sum_{j=q}^{n}|z_j|^2).
$$
Using (\ref{iEkappa}) we find then that
$$\langle i_E(\kappa),i_E(\kappa)\rangle
=\langle \omega_1\wedge\omega_2, \omega_1\wedge\omega_2\rangle
=\langle \omega_1,
\omega_1\rangle
\langle
\omega_2 , \omega_2
\rangle=(\sum_{j=1}^{q-1}|z_j|^2)
(\sum_{i=n-p+1}^{n}|z_i|^2).
$$
The integrals of the two terms in
(\ref{inn-a}) can be computed directly
and we prove then our results.
\end{proof}

\begin{theo+} Assume $p, q\ge 1$ and $p+q\le n$. 
\begin{enumerate}
\item
The  spectrum of the
 Laplacian 
$\square_{\mathcal Q} =\bar \partial_{\mathcal Q}^\ast 
\bar \partial_{\mathcal Q}
+\bar \partial_{\mathcal Q}\bar \partial_{\mathcal Q}^\ast 
$ on 
$L^2(\mathcal Q^{p, q-1})$ is given as follows:
 The 
 subspaces $ \Phi^{m, l}_{q-1, p}$  $(p+q <n)$,
$\Phi^{m, l}_{q-1, p+1}$, $\Phi^{m, l}_{q, p}$ ($p+q<n$),
and $\Phi^{m, l}_{q, p+1}$ $(p+q <n)$
are eigenspaces of
$\square_{\mathcal Q}$ with eigenvalues
$$
\frac{m+q-1}{m+n-p+1}
\left((l+n-q+1)(m+n-p)
-\frac1{n-1}(l+p)(m+q-1)\right),
$$
$$
(m+q-1)(l+n-q+2),\quad  (m+q-1)(l+n-q+1), \quad (m+q)(l+n-q+1)
$$
respectively.
\item
The spectra of the  Laplacians
$\triangle_{\mathcal Q} = 
\partial_{\mathcal Q}^\ast 
 \partial_{\mathcal Q}
+\partial_{\mathcal Q}\partial_{\mathcal Q}^\ast $
on 
$L^2(\mathcal Q^{p, q-1})$
are given
as follows: The 
 subspaces $ \Phi^{m, l}_{q-1, p}$  $(p+q <n)$,
$\Phi^{m, l}_{q-1, p+1}$ $(p+q <n)$, $\Phi^{m, l}_{q, p}$, 
and $\Phi^{m, l}_{q, p+1}$ $(p+q <n)$
are eigenspaces of
$\Delta_{\mathcal Q}$ with eigenvalues
$$
\frac{l+p}{l+n-q}
\left((l+n-q+1)(m+n-p)
-\frac1{n-1}(l+p)(m+q-1)\right),
$$
$$
(l+p)(m+n-p+1),\quad (l+p)(m+l-p+1), \quad (l+p)(m+n-p+1)
$$
respectively.
\end{enumerate}
\end{theo+}

\begin{proof} To simplify the notation we will
suppress the subindex $\mathcal Q$ and 
write $\partial=\partial_{\mathcal Q}$.
The operators $\partial$,
and $\bar \partial$, and their adjoint operators
are clearly $K$-invariant operators. By the multiplicity
free result in Proposition 2.4, we have
that the operator $\bar \partial_Q$
maps  highest weight vectors of bi-degree
$(p-1, q-1)$
to  highest weigh vectors of bi-degree
$(p, q-1)$, and the corresponding statement
 is true
 for the
operator $\partial_Q$ and
the adjoint operators. Explicitly
we have, by direct computations,
\begin{equation*}
\label{d-alpha}
\bar\partial \alpha=0, \quad
\partial\alpha=0
\end{equation*}

\begin{equation*}
\label{d-beta}
\bar\partial
\beta=0, \quad 
\partial
\beta_{q-1, p+1}=(l+p+1) \alpha_{q-1, p+1}, 
\end{equation*}

\begin{equation*}
\label{d-gamma}
\bar\partial \gamma_{q, p}=(m+q) \alpha_{q, p},
\quad 
\partial
\gamma_{q, p}=0
\end{equation*}
\begin{equation*}
\bar\partial
\delta_{q, p+1}^{m, l}
=(m+q)\beta_{q, p+1}^{m, l}, \quad
\partial
\delta_{q, p+1}^{m, l}
=(l+p+1)\gamma_{q, p+1}^{m, l}.
\end{equation*}
By the multiplicity free result we get then
$$
\bar\partial^\ast\alpha_{q-1, p}=(m+q-1)\frac{\Vert\alpha_{q-1, p}\Vert^2
}{\Vert\gamma_{q-1, p}\Vert^2
}\gamma_{q-1, p},\quad 
\partial^\ast\alpha_{q-1, p}=(l+p)\frac{\Vert
\alpha_{q-1, p}\Vert^2
}{\Vert\beta_{q-1, p}\Vert^2 }\beta_{q-1, p}
$$

$$
\bar\partial^\ast
\beta_{q-1, p+1}=(m+q-1)\frac{\Vert\beta_{q-1, p+1}\Vert^2}{\Vert \delta_{q-1, p+1}\Vert^2 }\delta_{q-1, p+1},\quad 
 \quad 
\partial^\ast \beta_{q-1, p+1}=0
$$

$$
\bar\partial^\ast
\gamma_{q, p}=0, \quad 
\partial^\ast \gamma_{q, p}=(l+p)
\frac{\Vert\gamma_{q, p}\Vert^2}{\Vert \delta_{q, p}\Vert^2 }\delta_{q, p},$$
$$
\bar\partial^\ast \delta_{q, p+1}=0,
\quad
\partial^\ast  \delta_{q-1, p+1}=0.
$$
 We compute
the eigenvalues of 
 $\square$ on $\alpha_{q-1, p}$, 
$$
\square\alpha_{q-1, p}
=\partial\bar\partial^\ast \alpha_{q-1, p}
=(m+q-1)\frac{\Vert\alpha_{q-1, p}\Vert^2
}{\Vert\gamma_{q-1, p}\Vert^2
}\partial
\gamma_{q-1, p}
=(m+q-1)^2\frac{\Vert\alpha_{q-1, p}\Vert^2
}{\Vert\gamma_{q-1, p}\Vert^2
}\alpha_{q-1, p},
$$
and our result then follows by Lemma 3.2.
The other eigenvalues can be computed similarly.
\end{proof}

\section{$SU(1, n)$-invariant
Unitary structure on  the harmonic
forms via the Szeg\"o map and the Tanaka Laplacians
}

In the paper \cite{Julg-Kasparov}
Julg and Kasparov give a geometric construction
of the $L^2$-harmonic $(p, q)$-forms on $B$. 
We recall some of their results here.

Adapting the notation there we put
$$
a=p+l, \quad b=q+m.
$$
Let $\omega=\bar z_1^l z_n^m d\bar z_1\wedge \cdots \wedge 
d\bar z_q\wedge dz_{n-q+1}\wedge \cdots \wedge dz_n $ be a $(p, q)$-form on $C^n$,
 which is
closed an primitive in the Euclidean metric
(thus Euclidean harmonic).  Let $F(t)={}_2F_1(a, b; a+b+2, t)$, $t=|z|^2$,
be the hypergeometric function. Put
$$
\widetilde \omega
= f_0\omega + f_1\tau\wedge i_T\omega +f_2\bar \tau\wedge
i_{\bar T}\omega + f_3\tau\wedge \bar \tau \wedge i_{\bar T}i_{T}\omega
$$
where
$$
f_0(z)=abF(t) -t F^\prime(t),
f_1(z)=(b+1)F^\prime(t), f_2(z)=(a+1)F^\prime(t), f_3(t)=F^{\prime\prime}(t).
$$

\begin{theo+}\cite[Sections 3-4]{Julg-Kasparov} 
Let $p, q\ge 1$ and $p+q=n$. The forms $\tilde\omega$ are $L^2$-harmonic
$(p, q)$-forms and they form an orthogonal basis
for the space $\mathcal H^{p, q}$. The Szeg\"o map
$S$ maps $L^2(\mathcal Q^{n-1})$ onto 
$\mathcal H^{p, q}$ and we have
$$
S(\gamma^{m, l}_{q, p})
=\frac{1}{aF(1)}\tilde \omega.
$$
\end{theo+}

In \cite[Theorem 4.9]{Julg-Kasparov} the unitarity
of the Szegö map is studied
on the space of all $L^2$-harmonic $n$-forms
using the Rumin operator, which is a fourth
order differential operator. We prove
now a unitarity result using the operator
$\partial_{\mathcal Q}$ and the vector field
$T+\bar T$. 

Define
$W^{1,2}_{p, q-1}$
to be the Hilbert space completion of differential
forms in $\Psi_{q, p}^{m,l}$, $m, l\ge 0$
with the norm
$$
 (\gamma, \gamma)_{W^{1,2}}
=( (\partial_{\mathcal Q}\partial_{\mathcal Q}
^\ast - (T+\bar T)_{\ast})\gamma, \gamma).
$$
Here  $(T+\bar T)_{\ast}$ is induced action of
 vector field $T+\bar T$ on
differential forms.
$W^{1,2}_{p, q-1}$ can be viewed as a Sobolev 
space of differential forms.

\begin{theo+}Let $p, q\ge 1$, $p+q=n$. The Szegö map
$S$ 
is a unitary operator from
$W^{1,2}_{p, q-1}$ onto the space $\mathcal H^{p, q}$ of
$L^2$-harmonic
$(p, q)$-forms.
\end{theo+}

\begin{proof} Using 
\cite{Julg-Kasparov}, Theorem 4.9 and Lemma 4.2 we can compute
easily the norm of $S(\gamma)$
in the $L^2(B, \wedge^{p, q})$-space of $(p, q)$-forms on the unit ball,
$$
(S(\gamma),S(\gamma))_{L^2(B, \wedge^{p, q})}
=\frac{m!l!}{(n)_{m+l+2}}b^2(a+1)(b+1),
$$
and by  Lemma 3.2 above
$$
(\gamma, \gamma)_{L^2(\mathcal Q^{p, q-1})}
=\frac{m!l!}{(n)_{m+l+2}}b(b+1).
$$
However
$$
-(T+\bar T)_{\ast} \gamma =  ((m+q)-(l+p))\gamma
= (b-a)\gamma
$$
by direct computation, and 
$$
\partial_{\mathcal Q}\partial_{\mathcal Q}
\gamma=a(b+1),
$$
by Theorem 3.3.
We have then
$$
(\partial_{\mathcal Q}\partial_{\mathcal Q}-(T+\bar T)_{\ast}) \gamma,
\gamma)=(S(\gamma), S(\gamma)),
$$
completing the proof.
\end{proof}

\begin{rema+} Theorem
4.2 needs to be modified slightly
for $(p, q)=(n, 0)$ or $(p, q)=(0, n)$. Consider the case
of $(p, q)=(n, 0)$. The space of harmonic
$(n, 0)$-forms are the 
 Bergman space
of holomorphic function $f(z)$ square integrable
with respect to the Lebesgue measure
with reproducing kernel $(1-\langle z, w\rangle)^{-(n+1)}$, viewed
  as $(n, 0)$-forms $f(z)(dz)^n$
where $(dz)^n$ stands for $dz_1\wedge \cdots \wedge dz_n$. 
Consider 
the space of $C^{\infty}$-functions 
on $S$  with the action of 
$G$:
$$
\pi(g): f(w)\mapsto f(g^{-1}w)J_{g^{-1}}(w)^{\frac{n}{n+1}}.
$$
The Szeg\"o map from functions on the sphere $S$
to the Bergman space is then
\begin{equation}
  \label{eq:szego-berg}
Sf(z)=\int_{S}\frac{1}{(1-\langle z, w\rangle)^{n+1}}
f(w)d\sig(w).
\end{equation}
The $L^2$-space on $S$ is decomposed
under $U(n)$ as a direct sum  of irreducible subspaces
with highest weights $(m, 0, \cdots, 0, -l)$, see
\cite{Rudin-ball}. Let
$f$ be in the subspace with weight $(0, \cdots, 0, -l)$,
 the Bergman space norm of $Sf$ is given
by
\begin{equation}
  \label{eq:iso-ber}
\Vert Sf\Vert^2(-\pi(T+\bar T)f, f)_{L^2(S)}  
\end{equation}
and thus extends to a unitary
operator from the corresponding
Sobolev space into the Bergman space.
Here
$T+\bar T$ is viewed as an element
in the complexification of the Lie algebra of $U(n)$
with the induced action.  This can be proved by a direct computation.
The above scalar  Szeg\"o map  is the same
as the  Szeg\"o map in Theorems 4.1 and 4.2; indeed
the area element is, writing
$$
d\sig(w)=\bar \tau\wedge (d\bar \tau)^n
=(d\bar w)^n \wedge (\sum_{j=1}^n (-1)^{j-1} w_j
d w_1\wedge \cdots \widehat{d w_j}\cdots \wedge d{w_n}),
$$
and identifying $Sf(z)$
in   (\ref{eq:szego-berg}) as $Sf(z)(dz)^n$
we can write   (\ref{eq:szego-berg})
 as
 \begin{equation*}
 \begin{split}
Sf(z)(dz)^n=\int_S \frac{(dz)^n
 \wedge
(d\bar w)^n}{(1-\langle z, w\rangle)^{(n+1)}}
 \wedge f(w) (\sum_{j=1}^n (-1)^{j-1} w_j
d w_1\wedge \cdots \widehat{d w_j}\cdots \wedge d{w_n}).   
\end{split}
 \end{equation*}
So if we identify the functions
$f(w)$ on the boundary with
$(n-1, 0)$-form
$F(w)=f(w)(\sum_{j=1}^n (-1)^{j-1} w_j
d w_1\wedge \cdots \widehat{d w_j}\cdots \wedge d{w_n})$; 
the scalar Szeg\"o map is then
a map on forms
$$S(F)(z)\int_{S}
 \frac{(dz)^n
 \wedge
(d\bar w)^n}{(1-\langle z, w\rangle)^{(n+1)}}
\wedge F(w),
$$
the integral kernel being then the $(n, n)$-form
$(1-\langle z, w\rangle)^{-(n+1)}(dz)^n
\wedge (d\bar w)^n$.
So   (\ref{eq:iso-ber}) takes
the form
\begin{equation}
  \label{eq:iso-ber-1}
\Vert SF\Vert^2(-(T+\bar T)_{\ast}F, F)_{L^2(S)}  
\end{equation}
with $(T+\bar T)_{\ast}$ the induced
action on differential forms
of the vector field.
\end{rema+}

\newcommand{\noopsort}[1]{} \newcommand{\printfirst}[2]{#1}
  \newcommand{\singleletter}[1]{#1} \newcommand{\switchargs}[2]{#2#1}
\providecommand{\bysame}{\leavevmode\hbox to3em{\hrulefill}\thinspace}
\providecommand{\MR}{\relax\ifhmode\unskip\space\fi MR }
\providecommand{\MRhref}[2]{%
  \href{http://www.ams.org/mathscinet-getitem?mr=#1}{#2}
}
\providecommand{\href}[2]{#2}

\end{document}